\documentclass{article}
\usepackage{cite}
\usepackage{amsmath,amsthm,amssymb}
\usepackage{verbatim}
\usepackage{graphicx}
\usepackage{footnote}

\newcommand{\RR}{\mathbb{R}} 
\newcommand{\mbf}{\mathbf}
\newcommand{\mtx}{\mbf}
\newcommand{\wh}{\widehat}

\title{The PowerURV algorithm for computing rank-revealing full factorizations}
\author{Abinand Gopal\footnote{Mathematical Institute, University of Oxford,
  email: \texttt{gopal@maths.ox.ac.uk}}~ and Per-Gunnar
  Martinsson\footnote{Institute for Computational Engineering and Sciences, UT
Austin, email: \texttt{pgm@ices.utexas.edu}}}
    \date{}
\begin{document}
\maketitle

\begin{center}
\begin{minipage}{120mm}
\textbf{Abstract:}
Many applications in scientific computing and data science require the
computation of a rank-revealing factorization of a large matrix. In many of
these instances the classical algorithms for computing the singular value
decomposition are prohibitively computationally expensive. The randomized
singular value decomposition can often be helpful, but is not effective unless
the numerical rank of the matrix is substantially smaller than the dimensions
of the matrix. We introduce a new randomized
algorithm for producing rank-revealing factorizations based on existing work 
by Demmel, Dumitriu and Holtz [\textit{Numerische Mathematik}, \textbf{108}(1), 2007]
that excels in this regime. The method is exceptionally easy to implement, and
results in close-to optimal low-rank approximations to a given matrix. The vast
majority of floating point operations are executed in level-3 BLAS, which leads to 
high computational speeds. The performance of the method is illustrated via several
numerical experiments that directly compare it to alternative techniques such as
the column pivoted QR factorization, or the QLP method by Stewart.
\end{minipage}
\end{center}

\section{Introduction} \label{s:intro}
In many applications such as in the analysis of large data sets and the
numerical solution of boundary integral equations, it is necessary to compute a
low-rank factorization of a matrix. One algorithm for this task that has proven
effective in the past several years is the randomized singular value
decomposition (RSVD) algorithm (see \cite{halko2011finding}).  Given a desired
truncation rank $k$, the RSVD computes a low-rank approximation by using a
random projection to reduce the dimensionality of the problem, paired with a
deterministic singular value decomposition (SVD) on the low-dimensional,
projected problem.  Often, a near-optimal low-rank approximation can be produced
in only a fraction of the time required to run a deterministic SVD on the entire
data matrix.  However, the RSVD has limitations. One is that it
requires information about the desired truncation rank prior to
computation, which may not be available. Another is that it becomes
uncompetitive when the target rank is not much smaller than the matrix dimensions. In these
regimes, which will be the focus of this report, it is usually helpful to
form a rank-revealing, full factorization of the matrix. By this we mean a
factorization of the original matrix that can easily be truncated to form
low-rank approximations.  The quintessential example of this is the full SVD,
which factors a matrix $\mbf{A} \in \RR^{m \times n}$ with $m \geq n$ into the
product $\mbf{A} = \mbf{U\Sigma V}^*$, where $\mbf{U} \in \RR^{m \times n}$ has
orthonormal columns, $\mbf{V} \in \RR^{n \times n}$ is orthogonal, and
$\mbf{\Sigma} \in \RR^{n \times n}$ is a diagonal matrix that takes values
$\mbf{\Sigma}(j,j) = \sigma_j$ for all $j = 1,2,\hdots,n$ such that $\sigma_1
\geq \sigma_2 \geq \hdots \sigma_n \geq 0$.  This factorization exists for any
matrix, and when it is obtained, a rank-$k$ approximation is given by $\mbf{A}_k
= \mbf{U}(:,1{:}k) \mbf{\Sigma}(1{:}k,1{:}k) \mbf{V}^*$.

The Eckart-Young theorem guarantees optimality in the spectral and
Frobenius norms of low-rank approximations obtained by truncating the SVD, but
computing an SVD can be prohibitively expensive in practice. A more economical
alternative is obtained through truncating a column pivoted QR factorization
(CPQR). While this is a much faster algorithm from both the perspective of
communication costs and operation count, there are no general guarantees on the
quality of the resulting low-rank approximation, and indeed there are known
cases where low-rank approximations obtained through truncating CPQR are
arbitrarily poor \cite{kahan1966numerical}.

A middle ground that can deliver results comparable in quality to the SVD, while
maintaining comparable efficiency to CPQR, can be achieved through the use of
so called UTV factorizations, which were introduced by G.W. Stewart in
\cite{stewart1992updating, stewart1993updating} and generalize the QR and SVD
factorizations. In full generality, the UTV factorization factors a given matrix
$\mbf{A} \in \RR^{m \times n}$ with $m \geq n$ into the product
\begin{equation}
  \mbf{A} = \mbf{U}\mbf{T}\mbf{V}^*,
  \label{eq:utv}
\end{equation}
where $\mbf{U} \in \RR^{m \times n}$ has orthonormal columns, $\mbf{V} \in
\RR^{n \times n}$ is orthogonal, and $\mbf{T} \in \RR^{n \times n}$ is
triangular.  While all of the algorithms we discuss in this report apply to both
the case where $\mbf{T}$ in \eqref{eq:utv} is upper-triangular and
lower-triangular, we will only consider the upper-triangular case for simplicity
and so will only consider URV decompositions. If we are given a matrix where $m
\leq n$, then we could do a ULV decomposition on $\mbf{A}^*$ and then transpose
the decomposition to find a URV decomposition. Thus, it also suffices to only
consider the case where $m \geq n$.

In order for a URV factorization to be useful for low-rank approximation, we
need the decomposition to be rank-revealing. By this we mean that for all $k =
1,2,\hdots,n$ a partition of the URV decomposition into
\begin{equation}
  \mbf{A} = \begin{bmatrix} \mbf{U}_1 & \mbf{U}_2 \end{bmatrix} \begin{bmatrix}
  \mbf{R}_{11} & \mbf{R}_{12} \\ & \mbf{R}_{22} \end{bmatrix} \begin{bmatrix}
  \mbf{V}_1^* \\ \mbf{V}_2^* \end{bmatrix},
  \label{eq:blocks}
\end{equation}
where $\mbf{U}_1 \in \RR^{m \times k}$, $\mbf{U}_2 \in \RR^{m \times
n-k}$, $\mbf{V}_1 \in \RR^{n \times k}$, $\mbf{V}_2 \in \RR^{n \times n
- k}$, $\mbf{R}_{11} \in \RR^{k \times k}$, $\mbf{R}_{12} \in \RR^{k
\times n - k}$ and $\mbf{R}_{22} \in \RR^{m-k \times n-k}$, has the
properties $\| \mbf{R}_{11} \|_2 \approx \sigma_k(\mbf{A})$ and $\| \mbf{R}_{22}
\|_2 \approx \sigma_{k+1}(\mbf{A})$.  This is not at all guaranteed for CPQR,
and although there do exist alternative pivoting strategies that guarantee
some type of rank-revealing property for QR (see \cite{gu1996efficient}), they are often
expensive to compute and still tend to deliver suboptimal results relative to
more general URV factorizations.

Given $\mbf{A}\in\RR^{m \times n}$ with $m \geq n$, a URV decomposition consists
of a matrix $\mbf{U} \in \RR^{m \times n}$ with orthonormal columns, an
orthogonal matrix $\mbf{V} \in \RR^{n \times n}$, and an upper-triangular matrix
$\mbf{R} \in \RR^{n \times n}$ such that $\mbf{A} = \mbf{URV}^*$. A CPQR
factorization of $\mbf{A} = \mbf{QRP}^{*}$ is a URV decomposition, where the
left orthogonal factor is given by $\mbf{Q}$ and the right orthogonal factor is
given by the permutation matrix $\mbf{P}$.  It follows that an unpivoted QR
factorization is also a URV factorization with the right orthogonal factor equal
to the $n \times n$ identity matrix.

\begin{figure}
  \centering
  \includegraphics[width=100mm]{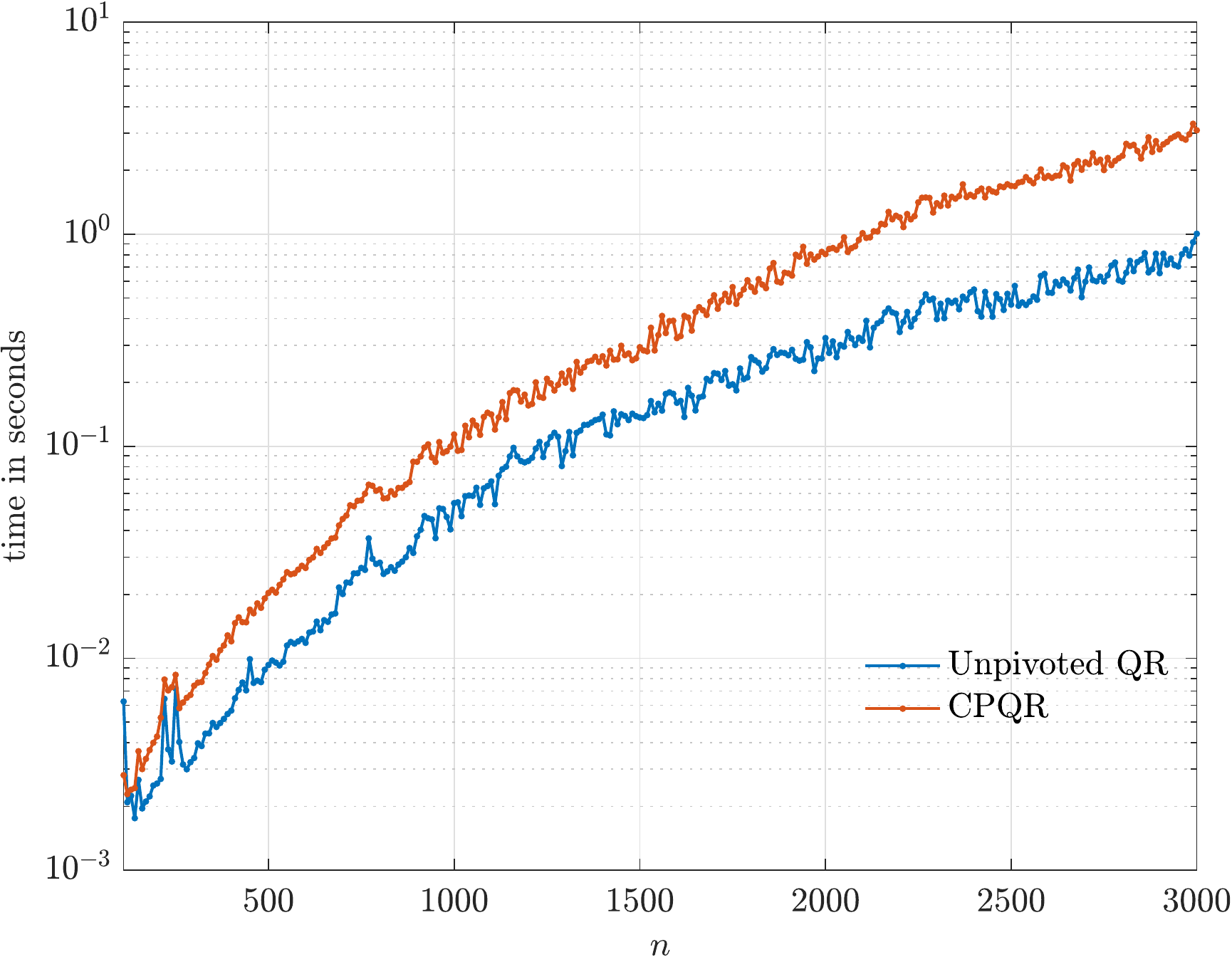}
  \caption{Comparison of wall clock time elapsed for unpivoted and column
  pivoted QR for Gaussian random matrices of size $n \times n$. Executed in
MATLAB on a 16-core Intel Xeon CPU E5-2643 @ 3.30 GHz with 64 GB of memory.}
  \label{fig:qrspeed}
\end{figure}

Throughout this work we will place an emphasis on communication costs.  While
low floating point operation counts are still crucial to efficient algorithms,
they do not tell the whole story on modern computing architectures.  This is
illustrated in a comparison between CPQR and unpivoted QR.  While both CPQR and
unpivoted QR have the same leading order floating point operation count, the
former is observed to be significantly slower than the latter as demonstrated in
Figure \ref{fig:qrspeed}.
This discrepancy is due to the fact that unpivoted QR can be executed almost
completely with level-3 BLAS operations, whereas CPQR requires a significant
number of level-2 BLAS operations to handle the pivoting.  We will be especially
interested in algorithms where the majority of floating point operations are in
matrix-matrix multiplications, which are intrinsically low in communication
costs. Moreover, matrix-matrix multiplication can easily leverage highly
parallel environments, such as GPUs, and can be further accelerated with fast
matrix-matrix algorithms, such as Strassen's algorithm.

\section{Algorithms for computing URV factorizations} \label{s:alg} As
aforementioned, the QR factorizations and the SVD are specific cases of URV
factorizations. We emphasize that we are primarily interested in URV
factorizations in the regime where low-rank approximations produced through CPQR
factorizations are too inaccurate and the SVD is too expensive to compute. Over
the past several decades many algorithms have been posed to compute the URV
factorization in this setting. In Section \ref{ss:stewart}, we summarize a
classical, deterministic algorithm for computing such a URV decomposition first
introduced by Stewart in \cite{stewart1999qlp}. In Section \ref{ss:demmel}, we
present another algorithm first introduced in the context of leveraging fast
matrix multiplication to accelerate eigenvalue computations by Demmel, Dumitriu,
and Holtz \cite{demmel2007fast} and then studied in its own right in
\cite{becker2017urv}. In Section \ref{ss:randutv}, we outline a recent,
randomized algorithm for computing the URV decomposition
\cite{martinsson2017randutv}.  In Section \ref{s:new}, we introduce a
modification of the algorithm by Demmel, Dumitriu, and Holtz which trades
slightly higher computational cost for stronger rank-revealing properties.

\subsection{Stewart's QLP} \label{ss:stewart}
One particularly effective algorithm for computing rank-revealing UTV
factorizations was introduced by Stewart in \cite{stewart1999qlp}, which
produces a ULV factorization through the use of two CPQR factorizations. Since
we have chosen to focus on the URV factorization in this report, we will instead
consider a variation where the QLP algorithm is applied to $\mbf{A}^*$ instead
of $\mbf{A}$. The transpose of the resulting factorization gives a URV
factorization for $\mbf{A}$.  To be precise, given $\mbf{A} \in \RR^{m
\times n}$, a CPQR factorization is performed on the transpose of $\mbf{A}$,
\[
  \mbf{A}^* = \mbf{Q}_1 \mbf{R}_1 \mbf{P}_1^*,
\]
where $\mbf{Q}_1 \in \RR^{n \times n}$ is an orthogonal matrix, $\mbf{R}_1 \in
\RR^{n \times m}$ is upper-trapezoidal, and $\mbf{P}_1 \in \RR^{m \times m}$ is
a permutation matrix. Then another CPQR factorization is
performed on the transpose of the product of the upper-trapezoidal factor and
the permutation matrix
\[
  (\mbf{R}_1 \mbf{P}_1^*)^* = \mbf{Q}_2 \mbf{R}_2 \mbf{P}_2^*,
\]
where $\mbf{Q}_2 \in \RR^{m \times n}$, $\mbf{R}_2 \in \RR^{n \times n}$ is
upper-triangular, and $\mbf{P}_2 \in \RR^{n \times n}$ is a permutation matrix.
Setting $\mbf{U} = \mbf{Q}_2$, $\mbf{R} = \mbf{R}_2$ and $\mbf{V} = \mbf{Q}_1
\mbf{P}_2$ yields a URV decomposition for $\mbf{A}$.

\begin{figure}[ht]
  \includegraphics[scale=0.4]{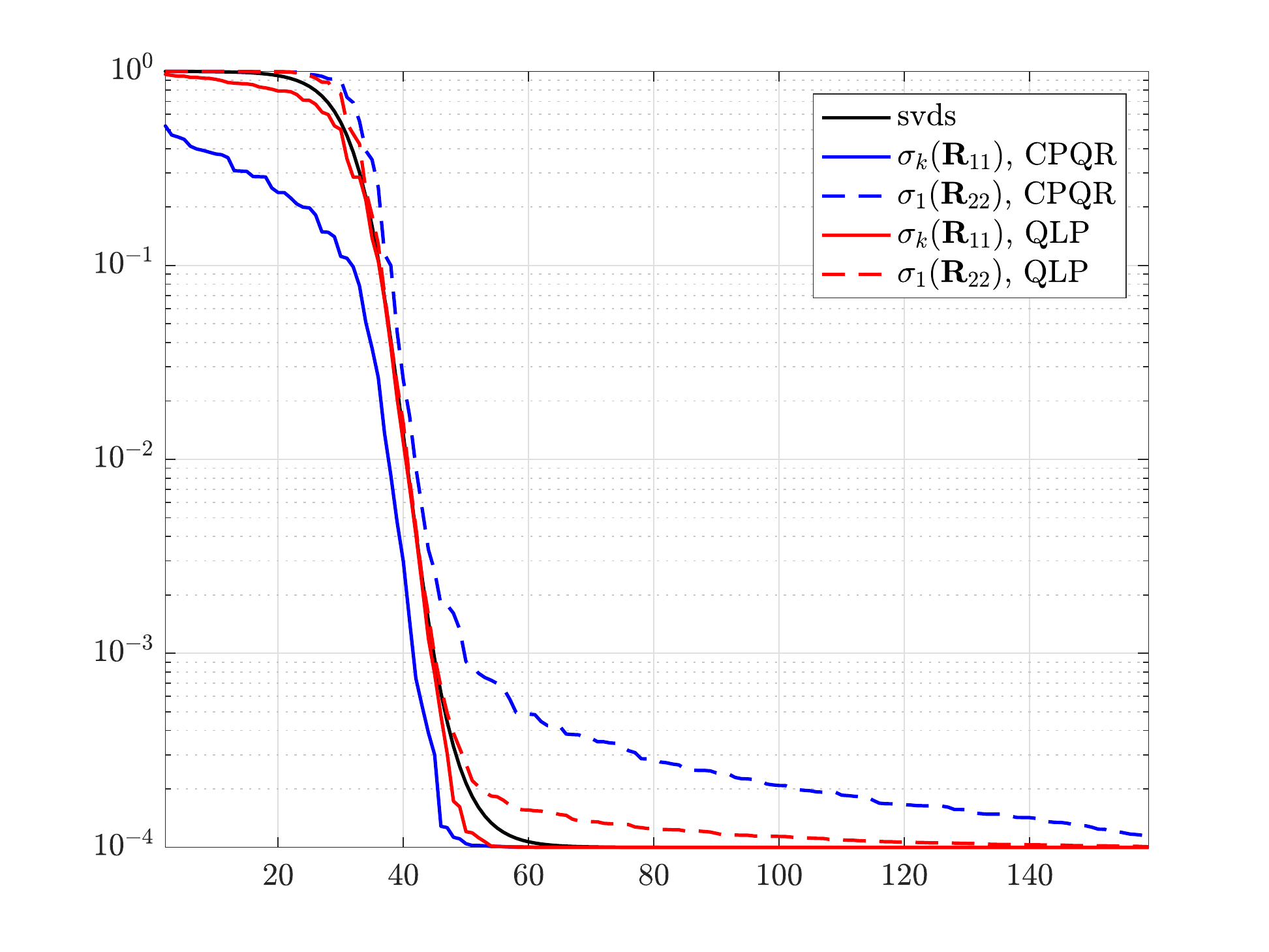} \hspace{-2em}
  \includegraphics[scale=0.4]{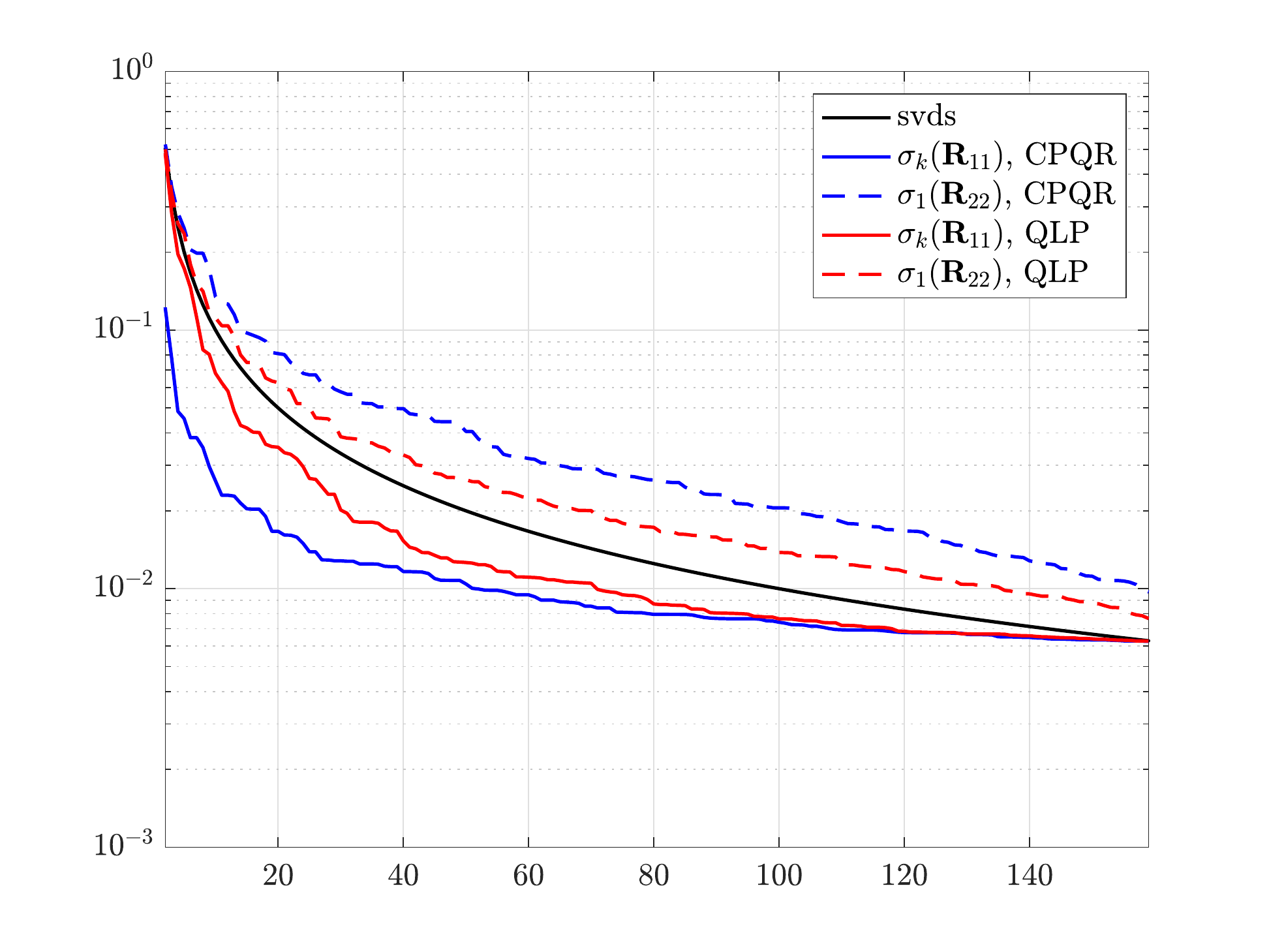}
  \caption{Comparison of qualities of rank-revealing factorizations obtained by
  CPQR and QLP for \emph{Matrix 3} (left) and \emph{Matrix 2} (right), both
specified in Section \ref{s:num}.}
  \label{fig:qlp}
\end{figure}

While the underlying mechanism behind this algorithm is the CPQR factorization,
it is empirically observed that using the second CPQR factorization produces a
significantly better rank-revealing decomposition. It was even shown that this
procedure could handle Kahan's example, for which CPQR by itself fails to
produce a rank-revealing factorization \cite{stewart1999qlp}. Figure
\ref{fig:qlp} illustrates the difference in quality of rank-revealing
factorizations for two non-pathological examples.

\subsection{DDH-URV} \label{ss:demmel}
The first randomized algorithm for computing the URV decomposition that we are
aware of was given in \cite{demmel2007fast} as a tool for accelerating
eigenvalue computations given a fast matrix-matrix multiplication algorithm. We
will refer to this algorithm in this technical report as the DDH-URV algorithm,
after the authors of \cite{demmel2007fast}. The algorithm starts by generating a
Haar distributed random matrix $\mbf{V}\in \RR^{n \times n}$ which is produced
by orthogonalizing a Gaussian random matrix using an unpivoted QR
factorization\footnote{Technically, a matrix produced in this fashion will
  usually not be exactly Haar distributed, since this would require restricting
  the upper-triangular factor in QR factorization to be non-negative along the
diagonal. In practice, this detail is unimportant, so we will continue to say
that orthogonal matrices produced through QR factorizations on a Gaussian random
matrix are Haar distributed. See \cite{edelman2005random} for more details.}.
Then, another unpivoted QR factorization is taken of the product $\mbf{AV}$ to
produce
\begin{equation}
  \mbf{AV} = \mbf{QR},
\label{eq:AV=QR}
\end{equation}
where $\mbf{Q} \in \RR^{m \times n}$ has orthonormal columns and $\mbf{R} \in
\RR^{n \times n}$ is upper-triangular.  We see that equation (\ref{eq:AV=QR})
directly yields the URV factorization $\mbf{A}=\mbf{QRV}^*$.

As in Stewart's QLP algorithm, two QR factorizations are required in this
algorithm. However, both QR factorizations are done through unpivoted QR, which
leads to a signficant saving in communication costs. This is because unpivoted
QR can operate on blocks, accumulate Householder reflectors, and then apply the
reflectors to the trailing matrix using a level-3 BLAS matrix-matrix multiply.
In contrast, CPQR needs a level-2 BLAS matrix-vector multiply
after orthonormalizing each column to update the trailing column norms and
choose the next pivot correctly. The price we pay for maintaining rank-revealing
properties while removing the pivoting is a matrix-matrix multiplication with a
dense Haar matrix, which can also be performed using a level-3 BLAS
matrix-matrix multiplication. This matrix-matrix multiplication can be further
accelerated by replacing the Haar matrix with a chain of random, structured
orthogonal matrices \cite{becker2017urv}.

While this algorithm is extremely computationally efficient, the quality of the
rank-revealing factorization produced can be rather low.  This is due to the
fact that the method does not use any information of the data matrix $\mbf{A}$
to form the approximation of the right singular space spanned by the columns of
$\mbf{V}$. Our numerical experiments show that low-rank approximations obtained
through truncation of a URV factorization produced with this algorithm are
substantially further from optimal than those produced than the other
algorithms discussed in this section.

\subsection{randUTV} \label{ss:randutv}
Another randomized algorithm for computing the URV decomposition in our regime
of interest is the randUTV algorithm presented in \cite{martinsson2017randutv}.
This algorithm is similar in structure to a blocked algorithm for computing a QR
factorization, but employs the idea of the RSVD to process the ``active'' block
of columns. The randUTV algorithm requires the selection of a power parameter $q
> 0$ and blocking parameter $b > 0$ a priori. It begins by drawing a Gaussian
random matrix $\mbf{G} \in \RR^{m \times b}$, followed by applying the power
iteration $\mbf{Y} = (\mbf{A}^* \mbf{A})^q \mbf{A}^* \mbf{G} \in \RR^{n \times
b}$, and then applying an unpivoted QR factorization on this result to obtain
the orthogonal matrix $\mbf{W} \in \RR^{n \times n}$. Then $\mbf{A}$ can be
written
\[
  \mbf{A} = \begin{bmatrix} \mbf{A} \mbf{W}_1 & \mbf{A} \mbf{W}_2 \end{bmatrix}
  \mbf{W}^*,
\]
where $\mbf{W}_1 = \mbf{W}(1:b,:) \in \RR^{n \times b}$ and $\mbf{W}_2 =
\mbf{W}(b+1:n,:) \in \RR^{n \times n-b}$. Then the SVD is taken
\[
  \mbf{A} \mbf{W}_1 = \wh{\mbf{U}} \mbf{D} \mbf{V}_s^*,
\]
where $\wh{\mbf{U}} \in \RR^{m \times b}$, $\mbf{D} \in \RR^{b \times b}$ and
$\mbf{V}_s \in \RR^{b \times b}$. We can then write
\begin{equation}
  \mbf{A} = \wh{\mbf{U}} \begin{bmatrix} \mbf{D} & \wh{\mbf{U}}^* \mbf{A}
  \mbf{W}_2 \end{bmatrix} \wh{\mbf{V}}^*,\quad \wh{\mbf{V}} = \mbf{W}
  \begin{bmatrix} \mbf{V}_s & \mbf{0} \\ \mbf{0} & \mbf{I}_{n-b} \end{bmatrix},
  \label{eq:randutv}
\end{equation}
where $\mbf{I}_{n-b}$ denotes the $(n-b) \times (n-b)$ identity matrix. We see
from \eqref{eq:randutv} that we have reduced the submatrix consisting of the
first $b$ columns of $\mbf{A}$ to an upper-trapezoidal matrix through left and
right multiplications by orthogonal matrices. We can then repeat the same
procedure on the bottom $(m-b) \times (n-b)$ submatrix of $\wh{\mbf{U}}^*
\mbf{A} \mbf{W}_2$ to obtain a factorization of $\mbf{A}$ where the submatrix
consisting of the first $2b$ columns are upper-trapezoidal and so on, until a
complete URV factorization is obtained. For more details, we refer the reader to
\cite{martinsson2017randutv}.

In practice we have found that randUTV can offer close to optimal low-rank
approximations at modest computational cost. With that said, it can be difficult
to implement. This is in contrast to the other algorithms described here, which
can all be implemented with just a few calls to standard BLAS routines. Another
subtlety of randUTV is that the rank-revealing quality is affected by the choice of
a blocking parameter, as opposed to just the computational time. Indeed, we have
observed that a poor choice of blocking parameter can produce poor
rank-revealing factorizations in addition to slowing down the computation.
Since we are not aware of any rigorous heuristics for determining this
parameter, we omit randUTV from our numerical experiments in Section
\ref{s:num}.

\section{PowerURV} \label{s:new}
The new PowerURV algorithm combines the power iteration in the RSVD with the
DDH-URV algorithm discussed in Section \ref{ss:demmel}. First a small integer $q
> 0$ is fixed, which controls how many steps of power iteration will be taken.
Then Gaussian random matrix $\mbf{G} \in \mathbb{R}^{n \times n}$ is drawn. We
then apply $q$ steps of a power iteration to $\mbf{G}$ and perform an unpivoted
QR factorization to produce
\begin{equation}
  \label{eq:power}
  (\mbf{A}^* \mbf{A})^q \mbf{G} = \mbf{VZ},
\end{equation}
where $\mbf{V} \in \RR^{n \times n}$ is orthogonal and $\mbf{Z} \in \RR^{n
\times n}$ is upper-triangular. An unpivoted QR factorization of the product
$\mbf{AV}$ is taken to obtain
\begin{equation}
  \mbf{AV}=\mbf{UR},
  \label{eq:powerqr}
\end{equation}
where $\mbf{U} \in \RR^{m \times n}$ has orthonormal columns and $\mbf{R} \in
\RR^{n \times n}$ is upper-triangular. A URV factorization is then given by
\[ \mbf{A} =\mbf{URV}^*. \]
In practice during the power iteration in \eqref{eq:power}, we need to
re-orthonormalize after each application of $\mbf{A}$ or $\mbf{A}^*$.
Otherwise, we see a loss of accuracy in the quality of our low-rank
approximations by $\varepsilon^{\frac{1}{2q+1}}$ where $\varepsilon$ denotes
machine precision. This is precisely analogous to the loss of accuracy in
subspace iteration when re-orthonormalization is omitted.

\begin{figure}[ht]
  \includegraphics[scale=0.4]{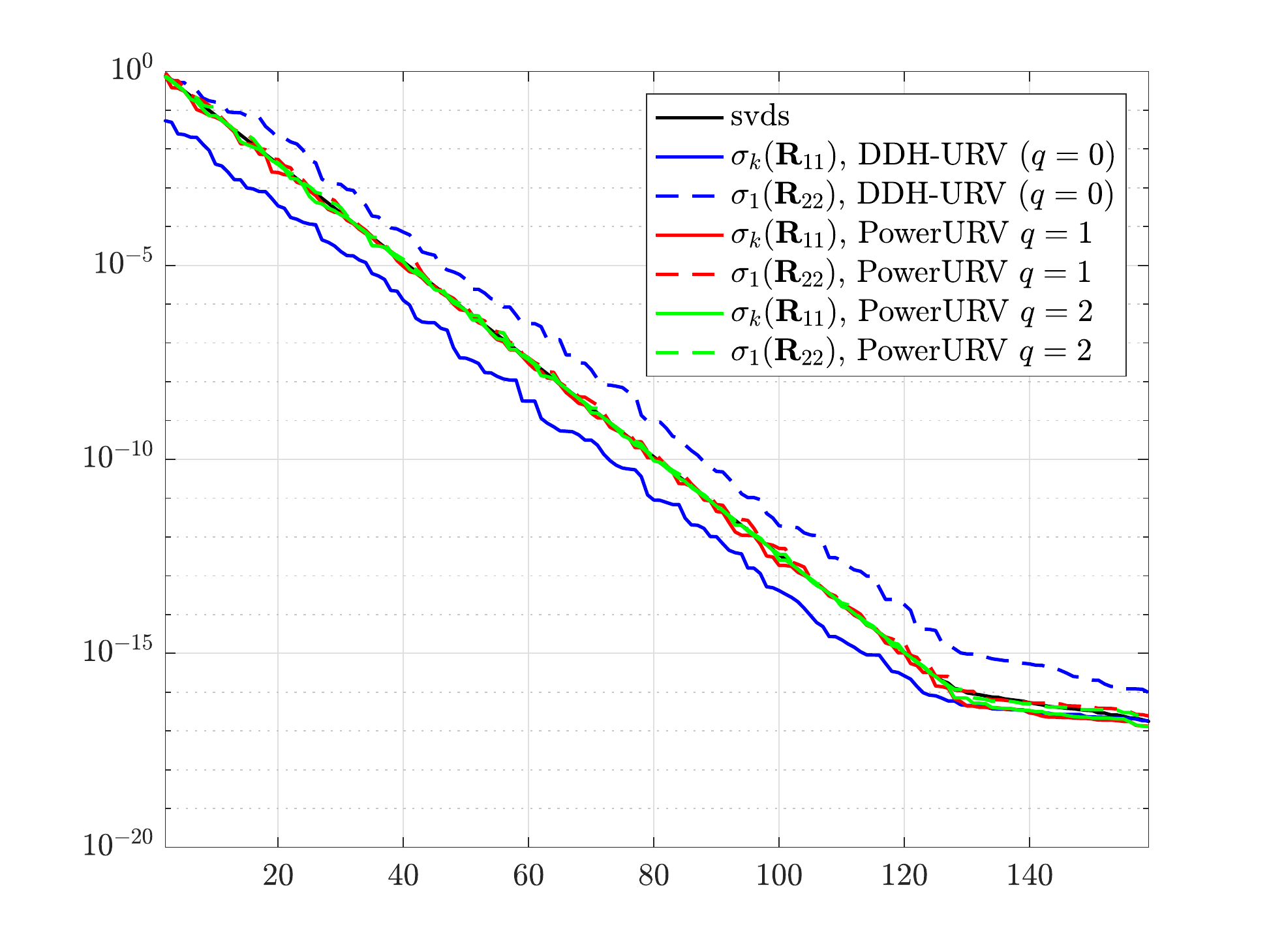} \hspace{-2em}
  \includegraphics[scale=0.4]{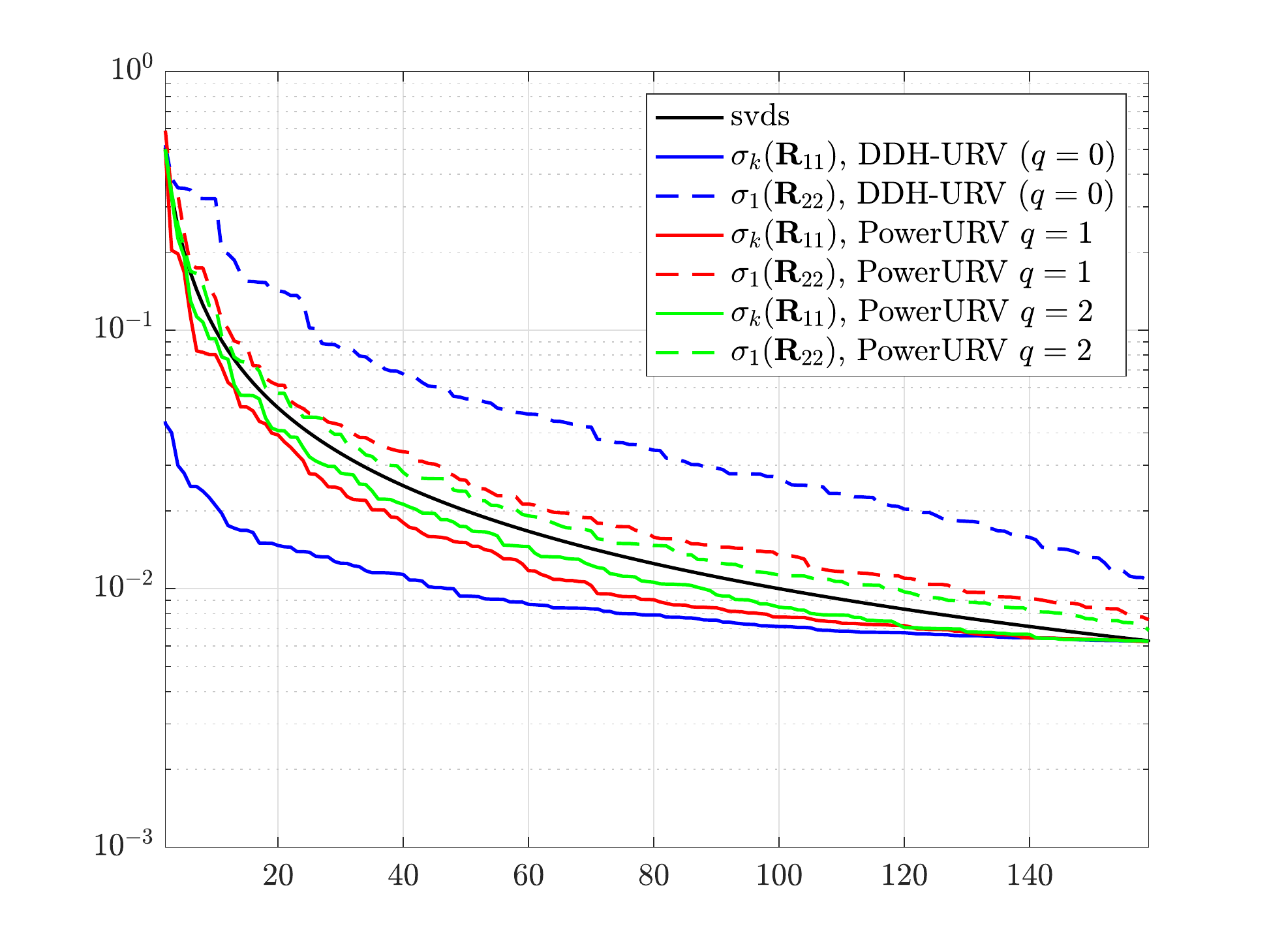}
  \label{fig:fastandslow}
  \caption{Comparison of qualities of rank-revealing factorizations obtained by
  PowerURV for $q=0,1,2$ for \emph{Matrix 1} (left) and \emph{Matrix 2} (right),
both specified in Section \ref{s:num}.}
\end{figure}

The PowerURV algorithm retains the same algorithmic structure as DDH-URV and,
when $q = 0$, is DDH-URV.  With $q > 0$, PowerURV takes into account information
of the row space of the data matrix $\mbf{A}$ when computing the approximation
to the space spanned by the right singular vectors. This allows for a
substantial improvement in quality when forming low-rank approximations. This
improvement is especially significant for matrices with slow decay in the
singular values. This is illustrated in Figure \ref{fig:fastandslow} which shows
the quantities relevant for rank-revealing for two 200$\times$160 matrices with
fast and slow singular value decay, respectively (see bullets on \emph{Matrix 1}
  and \emph{Matrix 2} in Section \ref{s:num} for the precise definitions of
these matrices). In the case of the matrix with fast decay, we see that the
$\sigma_1(\mbf{R}_{11})$ and $\sigma_k(\mbf{R}_{22})$ values of the
factorizations produced by the PowerURV algorithm with 1 and 2 power iterations
closely hug $\sigma_{k+1}(\mbf{A})$ for all $k=1,2,\hdots,n$. On the other hand,
the factorization without any of the power iterations is off by about an order
of magnitude in both directions. We see similar behavior for the matrix with
slow singular value decay except that there is a more discernible difference
between $q = 1$ and $q = 2$.

In contrast to Stewart's QLP algorithm, where the bulk of the operations
required are in two CPQR factorizations, the majority of work in PowerURV is in
level-3 BLAS matrix-matrix multiplications as in DDH-URV. This allows
PowerURV to be readily parallelized and reap advantages of modern computer
architectures. Another similarity to DDH-URV is that the implementation of
PowerURV is relatively straightforward with access to an efficient matrix-matrix
multiplication routine, an unpivoted QR factorization routine, and a random
number generator. We can defer choices of blocking parameters to the BLAS
routines themselves. The power iteration parameter $q$ depends only on spectral
information of the data matrix, and we find that in most cases it suffices to
set $q = 1$.

We also mention that the Gaussian matrix in \eqref{eq:power} could be replaced
by other random matrices such as a chain of orthogonal transforms (see
\cite{becker2017urv}).  However since this would come at a cost to accuracy and
only speed up 1 out of the $2q+1$ matrix-matrix multiplications, we prefer the
stronger theoretical guarantees of the Gaussian over the minute performance
advantages of using a single fast transform.

\section{Cost analysis} \label{s:cost}
\begin{figure}[ht]
  \centering
  \includegraphics[width=100mm]{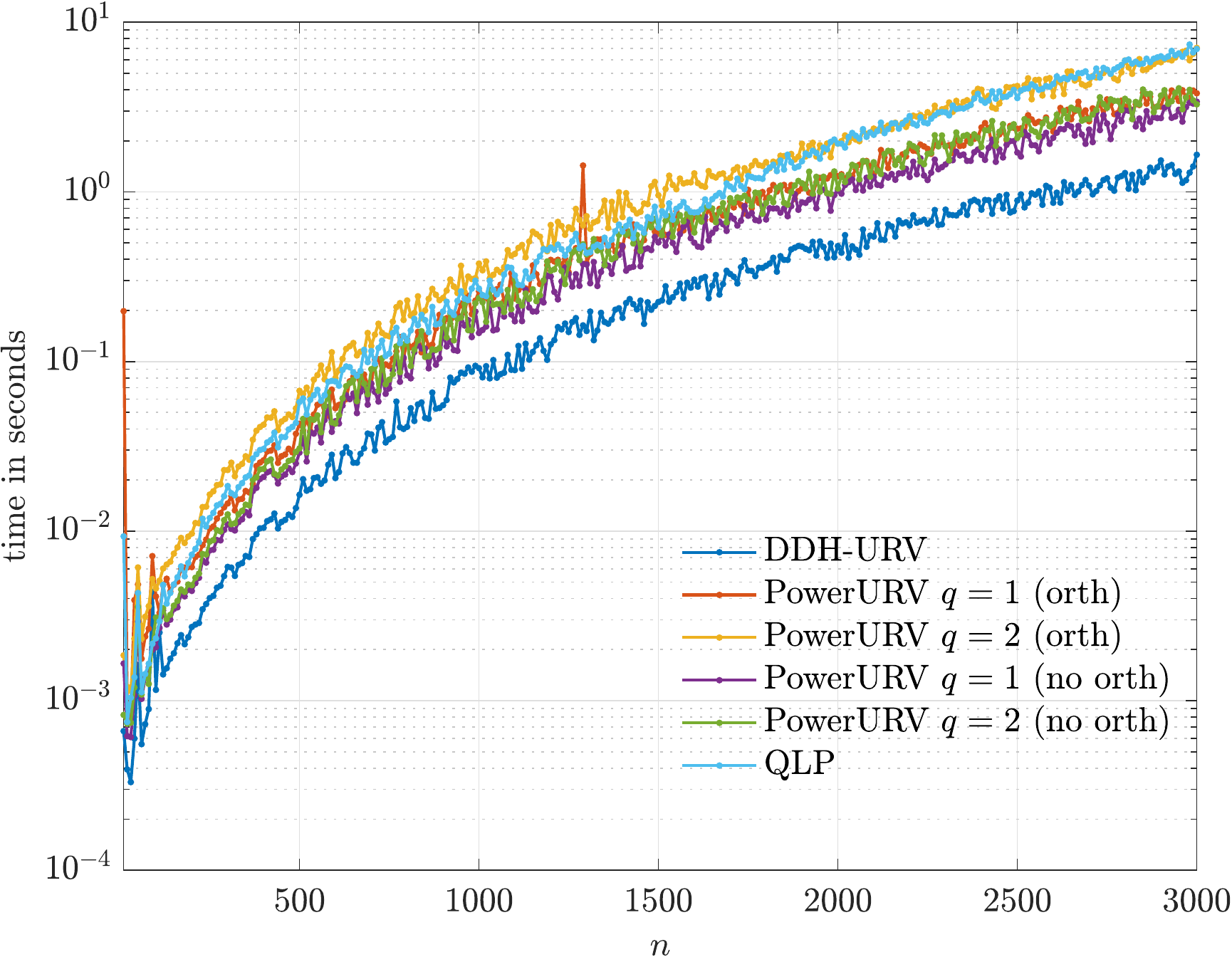}
  \caption{Comparison of wall clock time elapsed for the DDH-URV algorithm,
  PowerURV, and QLP for Gaussian random matrices of size $n \times n$. Executed
  in MATLAB on a 16-core Intel Xeon CPU E5-2643 @ 3.30 GHz with 64 GB of memory.}
  \label{fig:speedtest}
\end{figure}

Table \ref{t:cost} summarizes the leading order floating point operation counts
necessary for each of the algorithms discussed in this report for an $m \times
n$ matrix. For the SVD, we have only included the cost of bidiagonalization and
the costs of forming the left and right orthogonal factors. We note that the SVD
has a significantly higher operation count than all of the other algorithms
discussed.  Moreover, the bidiagonalization in the SVD cannot be done in level-3
BLAS. (This is a byproduct of having to apply Householder reflectors from both
the left and right one after the other). PowerURV has a significantly higher
flop count than the QLP algorithm, but as discussed earlier, flop counts alone
do not give the whole picture. The QLP algorithm needs two column pivoted QR
factorizations, which unlike the unpivoted QR factorization, needs a level-2
BLAS column norm update for correct pivot selection. As a result, despite the
higher flop count, we find in practice PowerURV to be as fast as QLP, if not
faster, for large problems, as illustrated in Figure \ref{fig:speedtest}.  This
also suggests that PowerURV has a better outlook for parallelism. The comparison
of the performance between randUTV and PowerURV is less clear. In our numerical
experiments, we found that the performance of randUTV is highly dependent on the
choice of blocking parameter, which can vary greatly based on hardware and
problem size.  It is clear that randUTV also has a high potential for
parallelism, but we believe that developing and tuning an efficient parallel
implementation would be non-trivial, whereas developing an effective parallel
implementation of PowerURV could consist of just calling parallel routines for
matrix-matrix multiplication and unpivoted QR factorization.

\begin{savenotes}
\begin{table}[ht]
  \centering
  \begin{tabular}{|c|c|c|c|}
    \hline
    Algorithm & DGEMM/ QR & CPQR & other level-2 BLAS\\ \hline
    Golub-Reinsch & 0 & 0 & $4m^2 n + 8mn^2 + 9n^3$  \\
    QLP & 0 & $2mn^2 + \frac{2}{3} n^3$ & 0 \\
    randUTV\footnote{The cost of forming the orthogonal factors has not been
    included.} & \small$(5+2q)mn^2 - \frac{1}{3}(3+2q)n^3$\normalsize
     & 0 & 0  \\
    PowerURV & \small$2(2q+1)m^2 n$\normalsize + \small$(4q+2)mn^2 - \frac{2}{3} (2q +
    1)n^3$\normalsize& 0 & 0 \\
    \hline
  \end{tabular}
  \label{t:cost}
\end{table}
\end{savenotes}

\section{Numerical results} \label{s:num}
\begin{figure}[ht]
  \centering
  \includegraphics[width=5.5in]{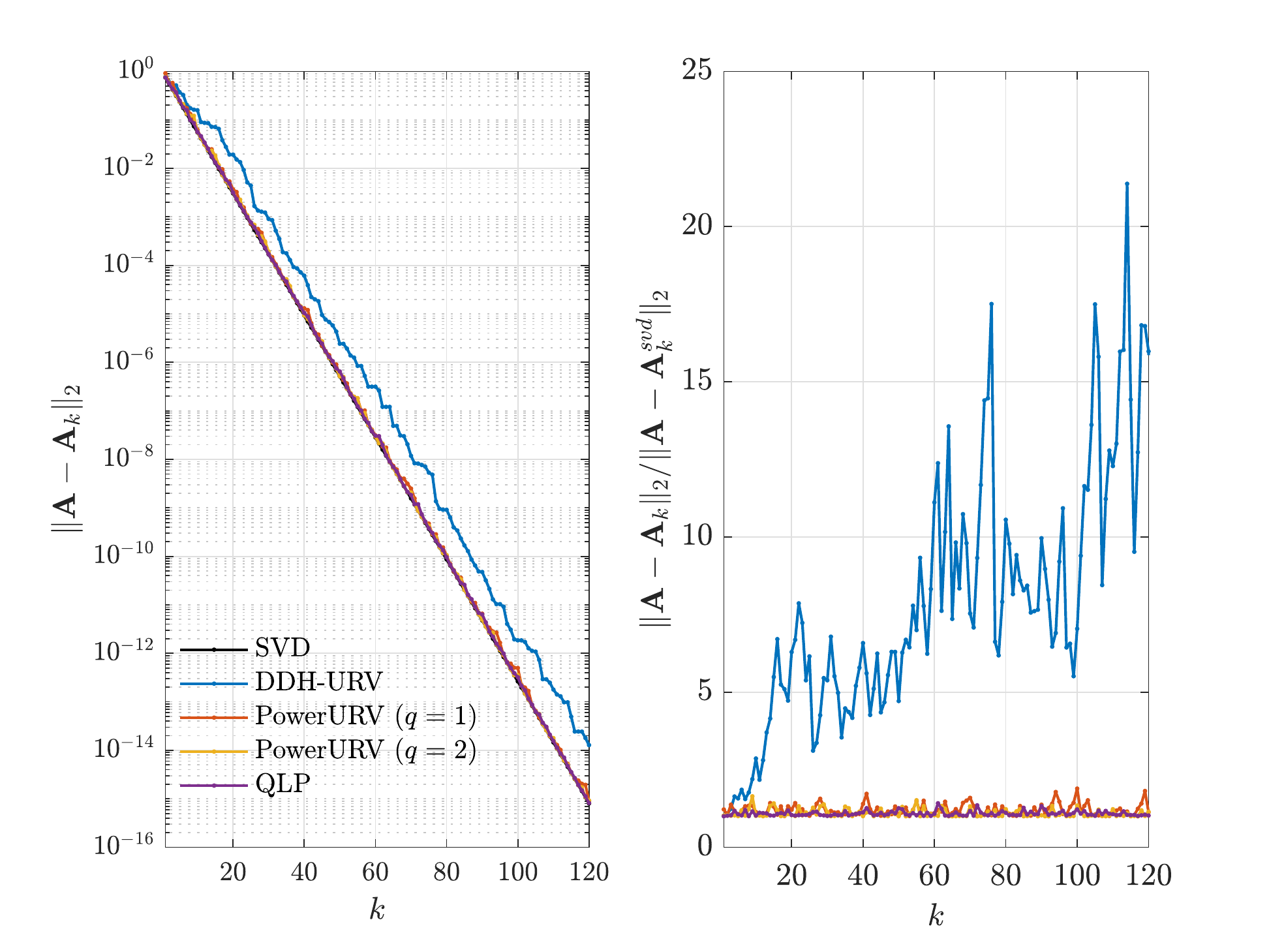}
  \caption{Rank-$k$ approximation errors to \emph{Matrix 1} in Section
  \ref{s:num}.}
  \label{fig:lr_fast}
\end{figure}

\begin{figure}
  \centering
  \includegraphics[width=5.5in]{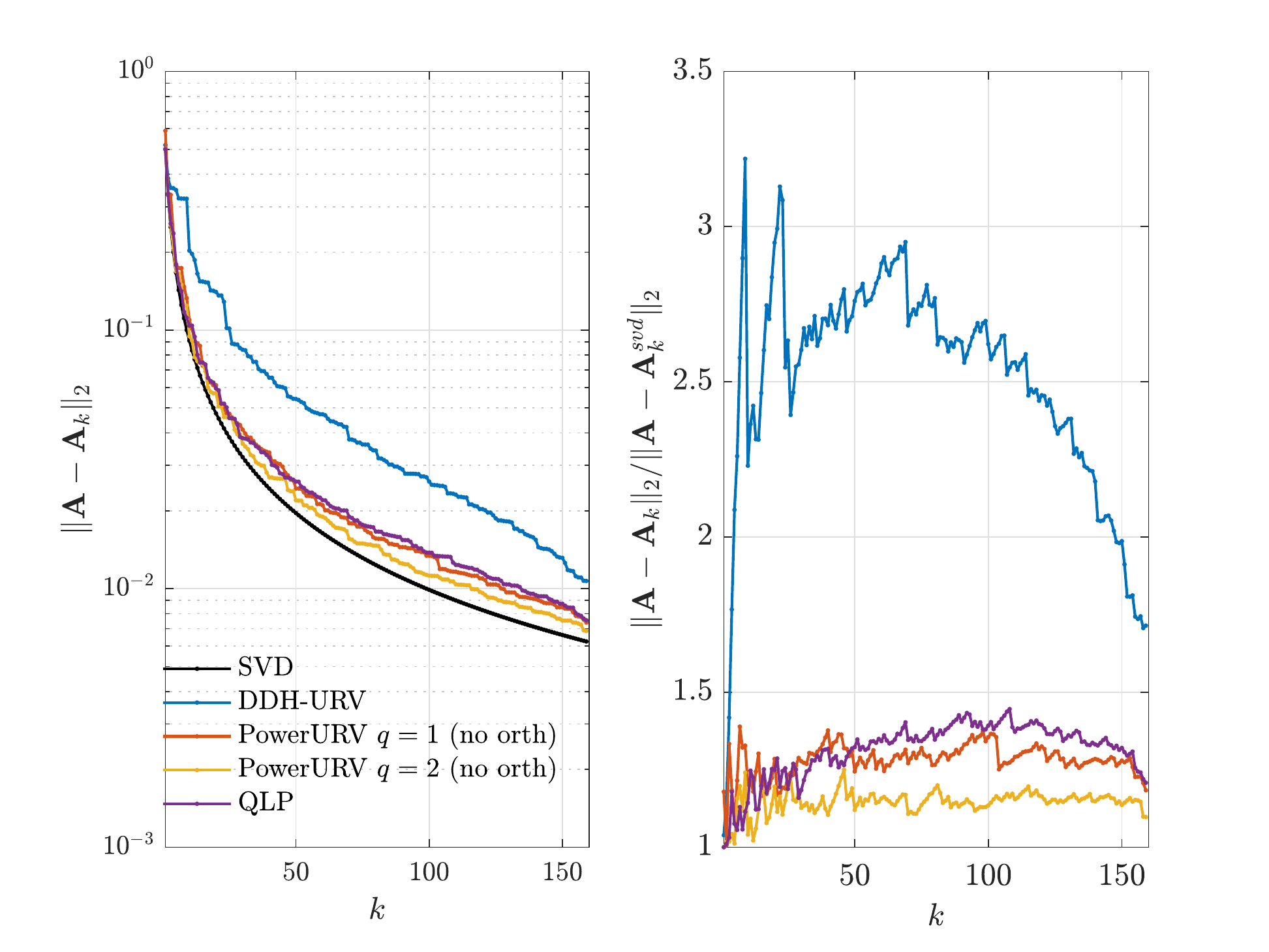}
  \caption{Rank-$k$ approximation errors to \emph{Matrix 2} in Section
  \ref{s:num}.}
  \label{fig:lr_slow}
\end{figure}

\begin{figure}
  \centering
  \includegraphics[width=5.5in]{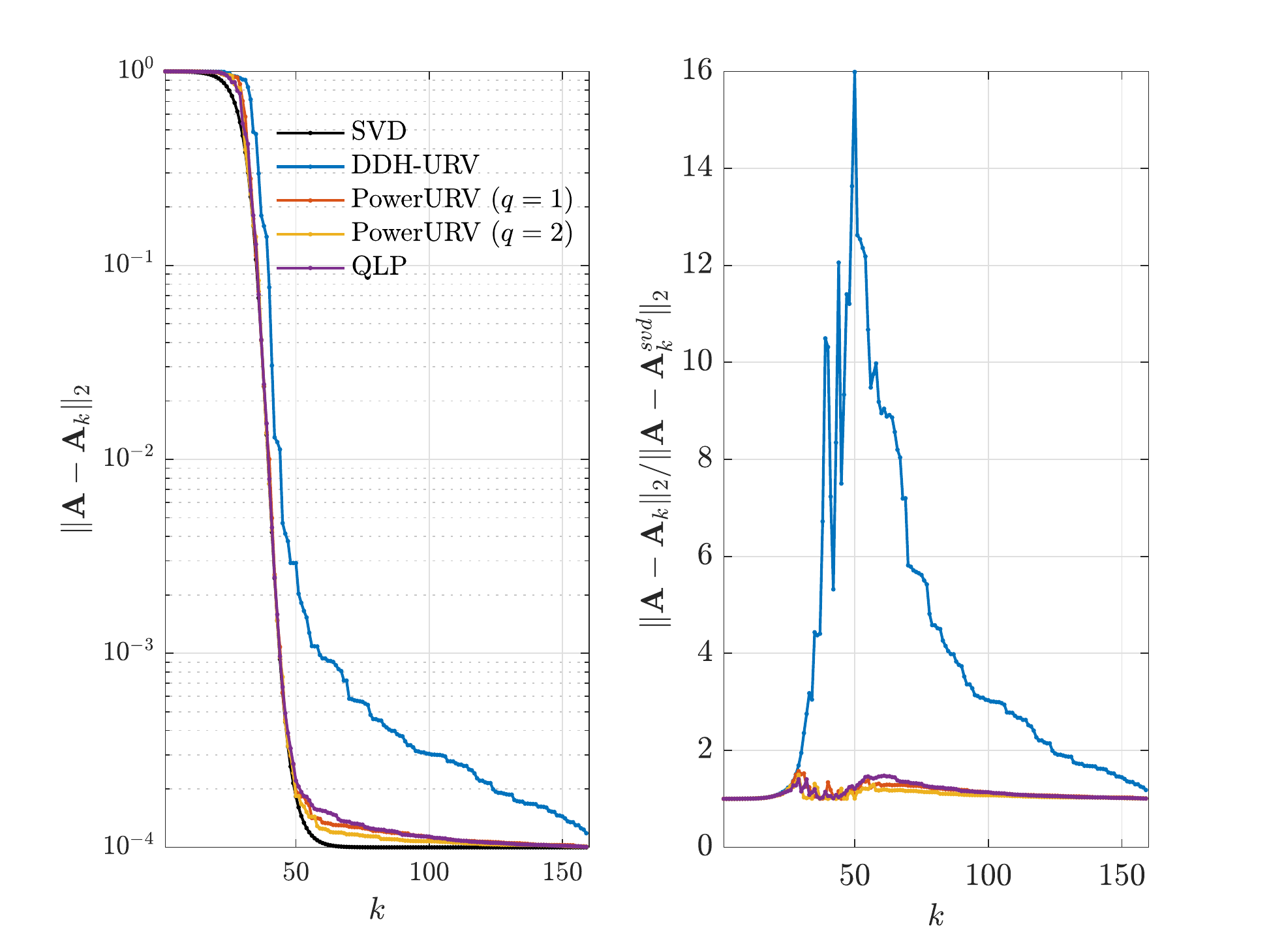}
  \caption{Rank-$k$ approximation errors to \emph{Matrix 3} in Section
  \ref{s:num}.}
  \label{fig:lr_sshape}
\end{figure}

\begin{figure}
  \centering
  \includegraphics[width=5.5in]{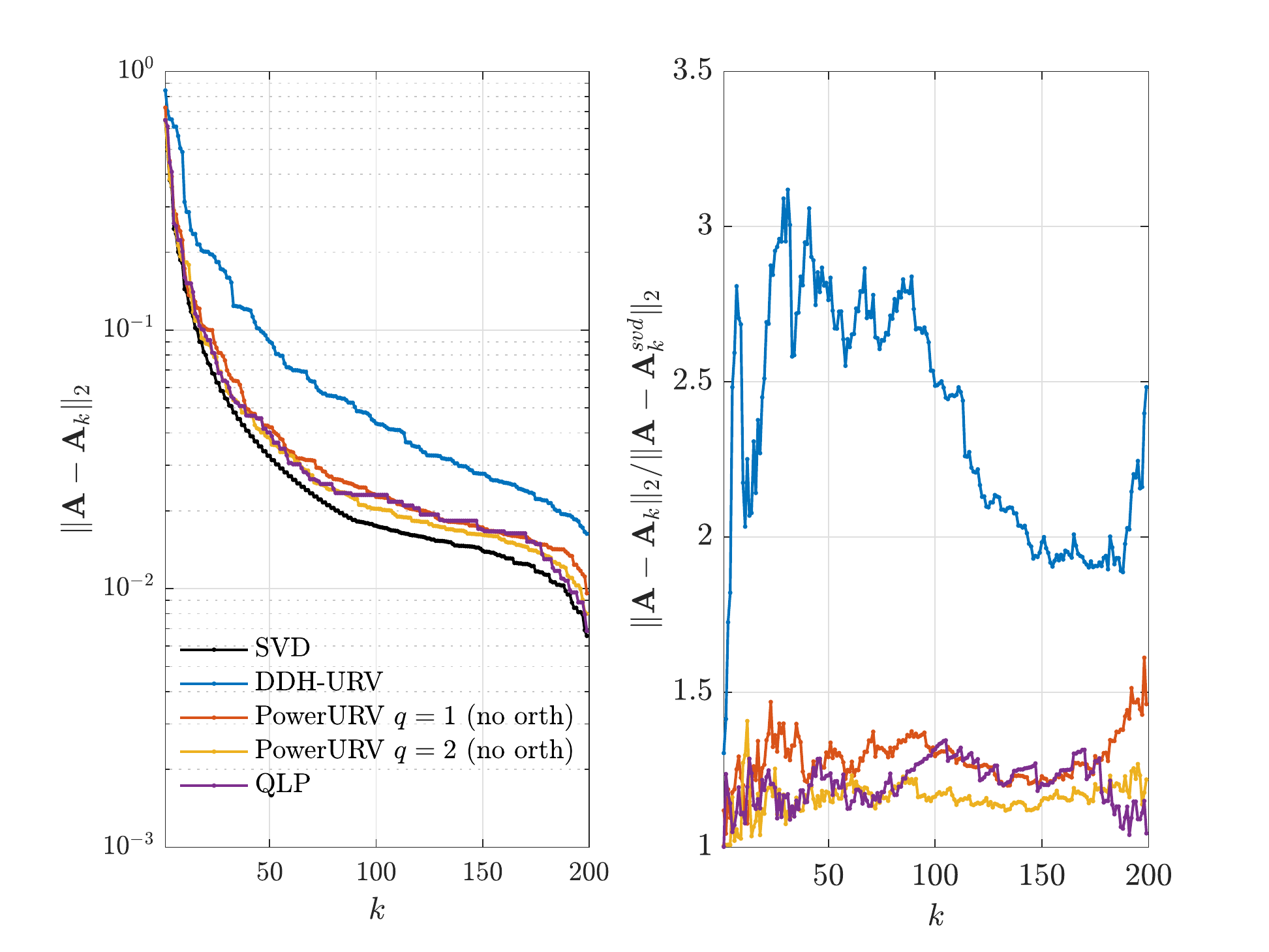}
  \caption{Rank-$k$ approximation errors to \emph{Matrix 4} in Section
  \ref{s:num}.}
  \label{fig:lr_bie}
\end{figure}

In this section, we compare the accuracy of QLP, DDH-URV, and PowerURV. We
follow the format in \cite{martinsson2017randutv} and benchmark the algorithms
against the following four matrices:
\begin{itemize}
    \item \emph{Matrix 1 (Fast Decay)}: This is a matrix $\mbf{A} = \mbf{UDV}^*
      \in \RR^{200 \times 160}$ where $\mbf{U}$ and $\mbf{V}$ are drawn from
      a Haar distribution and where $\mbf{D}$ is a rectangular, diagonal matrix
      with entries $\mbf{D}(k,k) = (10^{-20})^{k-1}$.
    \item \emph{Matrix 2 (Slow Decay)}: This is a matrix $\mbf{A} = \mbf{UDV}^*
      \in \RR^{200 \times 160}$ where $\mbf{U}$ and $\mbf{V}$ are drawn from
      a Haar distribution and where $\mbf{D}$ is a rectangular, diagonal matrix
      with entries $\mbf{D}(k,k) = k^{-1}$.
    \item \emph{Matrix 3 (S-Shaped Decay)}: This is a matrix $\mbf{A} =
      \mbf{UDV}^* \in \RR^{200 \times 160}$ where $\mbf{U}$ and $\mbf{V}$
      are drawn from a Haar distribution and where $\mbf{D}$ is a rectangular,
      diagonal matrix with entries $\mbf{D}(k,k) = 10^{-(1+\tanh(5(-1+2k/n)))}$
      for $k = 1,2,\hdots,80$ and $\mbf{D}(k,k) = 10^{-2}$ for $k =
      81,82,\hdots,150$.
    \item \emph{Matrix 4 (Boundary Integral Equation)}: This is a matrix
      $\mbf{A} \in \RR^{200 \times 200}$ that is the result of discretizing
      the boundary integral equation for the Laplace equation on a smooth,
      5-sided star.
\end{itemize}
In Figures \ref{fig:lr_fast}, \ref{fig:lr_slow}, \ref{fig:lr_sshape}, and
\ref{fig:lr_bie} have two plots for each matrix. For each URV algorithm, we plot
the relative and absolute errors of the difference between the $\mbf{A}$ and the
rank-$k$ approximation $\mbf{U}(:,1{:}k) \mbf{R}(1{:}k,:) \mbf{V}^*$. In all
cases we see that the PowerURV algorithm with just $q = 1$ compares very
favorably with the QLP algorithm. However, without any power iterations the
algorithm is a factor 2-10 worse. Since the singular values of \emph{Matrix 2}
and \emph{Matrix 4} do not decay below $\varepsilon^{1}{2q+1}$ for $q = 1,2$ we
omit the stabilizing orthogonalizations in the power iterations in these cases.

\section{Relationship with RSVD}
\label{sec:RSVDvsPURV}
The PowerURV algorithm is closely connected with the standard randomized
singular value decomposition algorithm (RSVD). To describe the connection,
let us briefly review the steps in the RSVD. Let $\mbf{A} \in \RR^{m \times n}$
with $m \geq n$. Given an integer $\ell < n$, the RSVD builds an
approximation to a truncated SVD via the following steps: A random matrix
(typically Gaussian) $\mbf{G}_{\rm rsvd} \in \RR^{n \times \ell}$ is drawn and the product
\begin{equation}
\label{eq:rsvd1}
  \mbf{Y}_{\rm rsvd} = \mbf{A}(\mbf{A}^* \mbf{A})^q \mbf{G}_{\rm rsvd} \in \RR^{m \times \ell},
\end{equation}
is evaluated. The columns $\mtx{Y}_{\rm rsvd}$ orthogonalized via an unpivoted QR factorization
\begin{equation}
\label{eq:rsvd2}
\mtx{Y}_{\rm rsvd} = \mtx{Q}_{\rm rsvd}\mtx{R}_{\rm rsvd}.
\end{equation}
In other words, the columns of the $m\times \ell$ matrix $\mtx{Q}_{\rm rsvd}$ form an orthogonal
basis for the  range of $\mtx{Y}_{\rm rsvd}$.
After this, a deterministic SVD of $\mbf{Q}_{\rm rsvd}^* \mbf{A}$ is computed to obtain
\begin{equation}
\label{eq:rsvd3}
\mbf{Q}_{\rm rsvd}^* \mbf{A} = \mbf{W}_{\rm rsvd}\mbf{\Sigma}_{\rm rsvd}(\mbf{V}_{\rm rsvd})^*,
\end{equation}
where $\mbf{W}_{\rm rsvd} \in \RR^{\ell \times \ell}$ is unitary,
where $\mbf{V}_{\rm rsvd} \in \RR^{n \times \ell}$ is orthogonal,
and $\mbf{\Sigma}_{\rm rsvd} \in \RR^{\ell \times \ell}$ is
diagonal with non-negative entries. The final step is to define the $m\times \ell$ matrix
\begin{equation}
\label{eq:rsvd4}
\mbf{U}_{\rm rsvd} = \mbf{Q}_{\rm rsvd} \mbf{W}_{\rm rsvd}.
\end{equation}
The end result is an approximate singular value decomposition
$$
\mbf{A} \approx \mbf{U}_{\rm rsvd}\mbf{\Sigma}_{\rm rsvd}\mbf{V}_{\rm rsvd}^{*}.
$$

The key claim in this section is that the first $\ell$ columns of the matrix $\mtx{U}$
resulting form PowerURV have exactly the same approximation accuracy as the columns of
the matrix $\mtx{U}_{\rm rsvd}$ resulting from the RSVD, provided that the same random
matrix is used. To be precise, we have:

\vspace{3mm}

\textit{\textbf{Lemma:}
Let $\mtx{A}$ be an $m\times n$ matrix, let $\ell$ be a positive integer such that $\ell < \min(m,n)$,
let $q$ be a positive integer,
and let $\mtx{G}$ be a matrix of size $n\times n$.
Let $\mtx{A} = \mtx{U}\mtx{R}\mtx{V}^{*}$ be the factorization resulting from the PowerURV algorithm,
as defined by (\ref{eq:power}) and (\ref{eq:powerqr}) and using $\mtx{G}$ as the starting point.
Let $\mtx{A} \approx \mtx{U}_{\rm rsvd}\mtx{\Sigma}_{\rm rsvd}\mtx{V}_{\rm rsvd}^{*}$ be the approximate
factorization resulting from RSVD, as defined
by (\ref{eq:rsvd1})-- (\ref{eq:rsvd4})), starting with $\mtx{G}_{\rm rsvd} = \mtx{G}(:,1:\ell)$.
Suppose that the rank of $\mtx{A}(\mtx{A}^{*}\mtx{A})^{q}\mtx{G}_{\rm rsvd}$ is no lower than the
rank of $\mtx{A}$
(this holds with probability 1 when $\mtx{G}_{\rm rsvd}$ is Gaussian).
Then
$$
\mtx{U}(:,1:\ell)\mtx{U}(:,1:\ell)^{*}\mtx{A} =
\mtx{U}_{\rm rsvd}\mtx{U}_{\rm rsvd}^{*}\mtx{A}.
$$
}

\begin{proof} We can without loss of accuracy assume that the matrix
$\mtx{A}$ has rank at least $\ell$.
(If it is rank deficient, then the proof we give will apply for a modified
$\ell' = \mbox{rank}(\mtx{A})$, and it is easy to see that adding additional
columns to the basis matrices will make no difference since in this case
$\mtx{U}(:,1:\ell)\mtx{U}(:,1:\ell)^{*}\mtx{A} = \mtx{U}_{\rm rsvd}\mtx{U}_{\rm rsvd}^{*}\mtx{A} = \mtx{A}$.)

We will prove that $\mbox{Ran}(\mtx{U}(:,\ell)) = \mbox{Ran}(\mtx{U}_{\rm rsvd})$, which
immediately implies that the projectors
$\mtx{U}(:,1:\ell)\mtx{U}(:,1:\ell)^{*}$ and
$\mtx{U}_{\rm rsvd}\mtx{U}_{\rm rsvd}^{*}$ are identical. Let us first
observe that restricting (\ref{eq:power}) to the first $\ell$ columns, we obtain
\begin{equation}
\label{eq:swim0}
(\mtx{A}^{*}\mtx{A})^{q}\mtx{G}_{\rm rsvd} =
\mtx{V}(:,1:\ell)\mtx{Z}(1:\ell,1:\ell).
\end{equation}
We can then connect $\mtx{Y}_{\rm rsvd}$ and $\mtx{U}(:,1:\ell)$ via a simple computation
\begin{equation}
\label{eq:swim1}
\mtx{Y}_{\rm rsvd} \stackrel{(\ref{eq:rsvd1})}{=}
\mtx{A}(\mtx{A}^{*}\mtx{A})^{q}\mtx{G}_{\rm rsvd} \stackrel{(\ref{eq:swim0})}{=}
\mtx{A}\mtx{V}(:,1:\ell)\mtx{Z}(1:\ell,1:\ell) \stackrel{(\ref{eq:powerqr})}{=}
\mtx{U}(:,1:\ell)\mtx{R}(1:\ell,1:\ell)\mtx{Z}(1:\ell,1:\ell).
\end{equation}
Next we link $\mtx{U}_{\rm rsvd}$ and $\mtx{Y}_{\rm rsvd}$ via
\begin{equation}
\label{eq:swim2}
\mtx{U}_{\rm rsvd} \stackrel{(\ref{eq:rsvd4})}{=}
\mtx{Q}_{\rm rsvd}\mtx{W}_{\rm rsvd} \stackrel{(\ref{eq:rsvd2})}{=}
\mtx{Y}_{\rm rsvd}\mtx{R}_{\rm rsvd}^{-1}\mtx{W}_{\rm rsvd}.
\end{equation}
Combining (\ref{eq:swim1}) and (\ref{eq:swim2}), we find that
\begin{equation}
\label{eq:swim3}
\mtx{U}_{\rm rsvd} =
\mtx{U}(:,1:\ell)\mtx{R}(1:\ell,1:\ell)\mtx{Z}(1:\ell,1:\ell)
\mtx{R}_{\rm rsvd}^{-1}\mtx{W}_{\rm rsvd}.
\end{equation}
The rank assumption implies that the $\ell\times \ell$ matrix
$\mtx{R}(1:\ell,1:\ell)\mtx{Z}(1:\ell,1:\ell)\mtx{R}_{\rm rsvd}^{-1}\mtx{W}_{\rm rsvd}$
is non-singular, which establishes that the matrices $\mtx{U}_{\rm rsvd}$ and
$\mtx{U}(:,1:\ell)$ have the same range.
\end{proof}

\begin{figure}
\centering
\includegraphics[width=120mm]{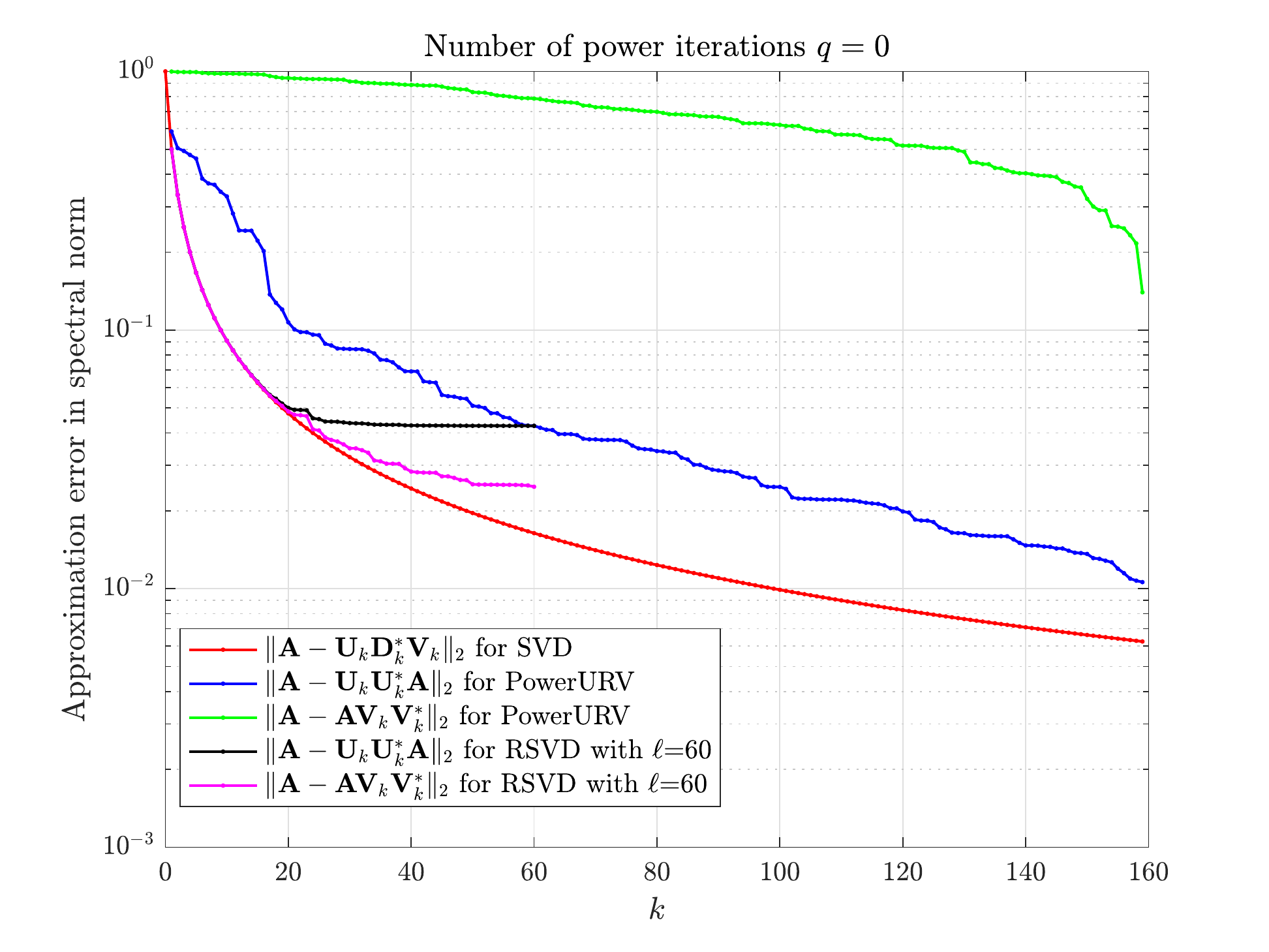}
\caption{Comparison of how well the factorizations resulting from the
RSVD and PowerURV algorithms reveal numerical rank, as discussed
in Section \ref{sec:RSVDvsPURV}. The matrix $\mtx{A}$ is ``Matrix 2''
described in Section \ref{s:num}, and the RSVD was executed with $\ell=60$.
No power iteration was used for either method (so $q=0$).}
\label{fig:compare_q=0}
\end{figure}

\begin{figure}
\centering
\includegraphics[width=120mm]{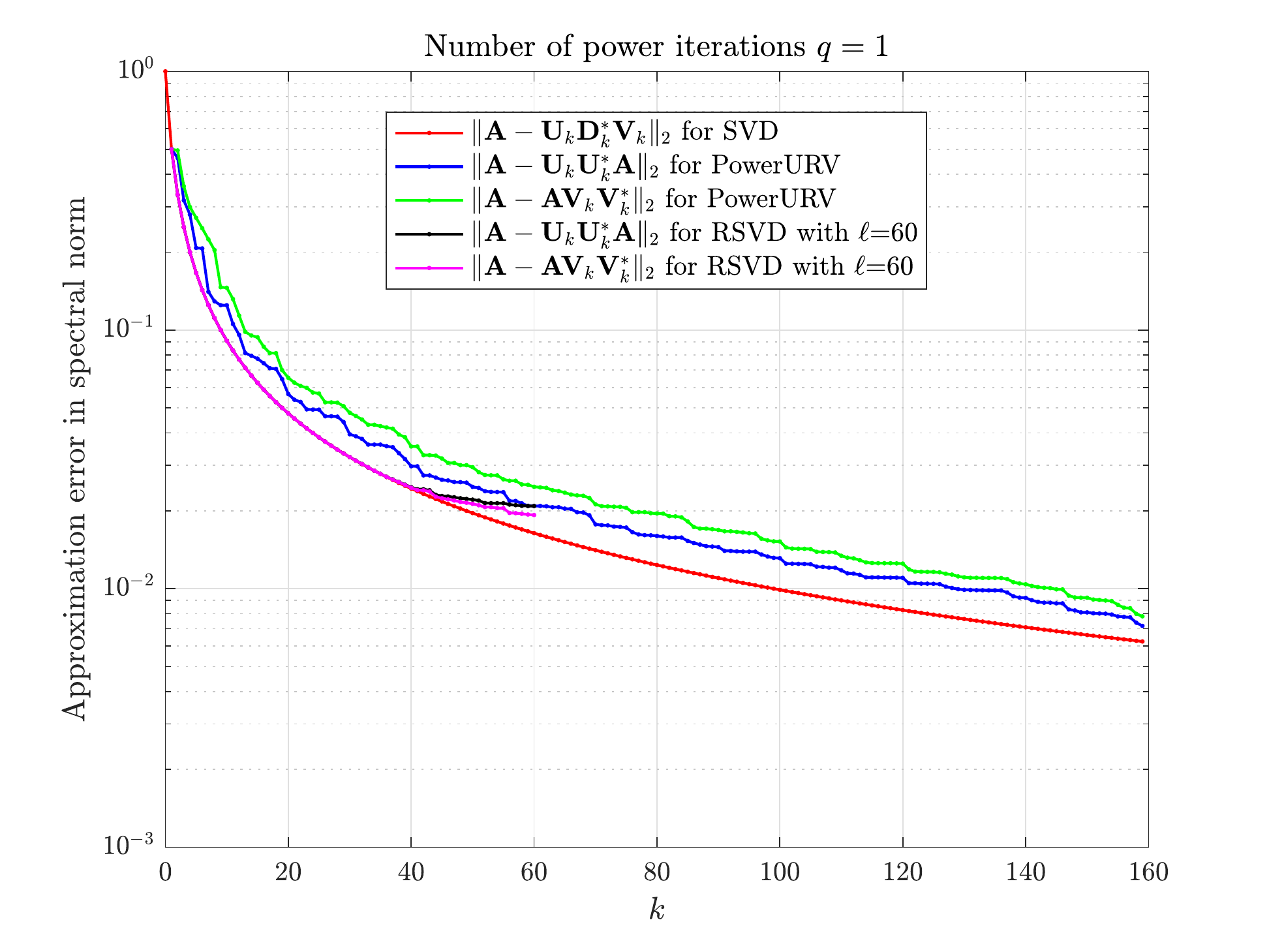}
\caption{The same experiment as shown in Figure \ref{fig:compare_q=0}, but
now with one step of power iteration, so that $q=1$.}
\label{fig:compare_q=1}
\end{figure}

The equivalency established in the Lemma between RSVD and PowerURV allows for much of the
theory for analyzing the RSVD in \cite{halko2011finding} to directly apply to the PowerURV
algorithm. To illustrate the theorem with a specific example, we show in Figures
\ref{fig:compare_q=0} and \ref{fig:compare_q=1} how well the columns of
$\mtx{U}$ and $\mtx{U}_{\rm rsvd}$ span the column space of the ``Matrix 2 (slow decay)''
we introduced in Section \ref{s:num}. To be precise, we ran the RSVD with $\ell = 60$,
and plotted the approximation errors
$$
\|\mtx{A} - \mtx{U}(:,1:k)\mtx{U}(:,1:k)^{*}\mtx{A}\|
\qquad\mbox{and}\qquad
\|\mtx{A} - \mtx{U}_{\rm rsvd}(:,1:k)\mtx{U}_{\rm rsvd}(:,1:k)^{*}\mtx{A}\|
$$
as a function of $k$, for $k \in \{1,2,\dots,\ell\}$. The equivalence claimed
in the lemma is captured by intersections of the errors at $k=\ell=60$.

It is important to note that while there is a close connection between RSVD and PowerURV,
the RSVD is able to attain substantially higher overall accuracy than PowerURV since it
can take advantage of one additional application of $\mtx{A}$. (To wit, RSVD requires
$2q+2$ applications of either $\mtx{A}$ or $\mtx{A}^{*}$, while PowerURV requires only $2q+1$.)
This additional application makes the first $k$ columns of $\mtx{U}_{\rm rsvd}$ a much
better basis than the first $k$ columns of $\mtx{U}$, as long as $k$ does not get close to $\ell$.
(One might say that the matrix $\mtx{W}_{\rm rsvd}$ rearranges the columns inside $\mtx{Q}_{\rm rsvd}$
to make the leading columns much better aligned with the corresponding singular vectors.)
This effect is visible in Figures \ref{fig:compare_q=0} and \ref{fig:compare_q=1} by the fact
that the black dots are initially much closer to the minimal errors on the red line than the
blue dots. It is only at the very end that the two lines meet.
An additional way that the RSVD benefits from the additional application of $\mtx{A}$ is that
the columns of $\mtx{V}_{\rm rsvd}$ end up being a far more accurate basis for the row space of
$\mtx{A}$ than the columns of the matrix $\mtx{V}$ resulting from the PowerURV. This was of course
expected since for $q=0$, the matrix $\mtx{V}$ incorporates no information from $\mtx{A}$ at all.
In Figures \ref{fig:compare_q=0} and \ref{fig:compare_q=1} we see this effect by noticing how much
smaller the errors marked by the magenta lines are than the errors marked by the green lines.

\section{Conclusion} \label{s:conc}

We described the randomized algorithm PowerURV for computing a rank-revealing factorization
of a general matrix. The method is computationally efficient since it relies only
on matrix-matrix multiplications and unpivoted QR factorizations involving the full
matrix. This makes it highly efficient on modern communication constrained hardware.
It is also exceptionally simple to implement. The method builds off of existing work
by Demmel, Dumtriu, and Holtz \cite{demmel2007fast}, and also ties in to work on
the Randomized SVD \cite{halko2011finding}. We compared the speed and accuracy
of PowerURV to previously proposed algorithms for computing rank-revealing factorizations,
such as the full SVD, the so called QLP method by Stewart \cite{stewart1999qlp}, and the
original randomized method proposed by Demmel, Dumtriu, and Holtz. In the regime where
a full factorization is sought (as opposed to cases where the numerical rank of the matrix
is far smaller than the matrix dimensions), PowerURV provides an excellent compromise between
computational speed on the one hand, and quality in terms of the rank-revealing properties
on the other.

\pagebreak
\bibliographystyle{abbrv}
\bibliography{gopal}

\begin{thebibliography}{10}

\bibitem{becker2017urv}
S.~Becker, J.~Folberth, and L.~Grigori.
\newblock {URV} factorization with random orthogonal system mixing.
\newblock {\em arXiv preprint arXiv:1703.02499}, 2017.

\bibitem{demmel2007fast}
J.~Demmel, I.~Dumitriu, and O.~Holtz.
\newblock Fast linear algebra is stable.
\newblock {\em Numerische Mathematik}, 108(1):59--91, 2007.

\bibitem{edelman2005random}
A.~Edelman and N.~R. Rao.
\newblock Random matrix theory.
\newblock {\em Acta Numerica}, 14:233--297, 2005.

\bibitem{gu1996efficient}
M.~Gu and S.~C. Eisenstat.
\newblock Efficient algorithms for computing a strong rank-revealing {QR}
  factorization.
\newblock {\em SIAM Journal on Scientific Computing}, 17(4):848--869, 1996.

\bibitem{halko2011finding}
N.~Halko, P.-G. Martinsson, and J.~A. Tropp.
\newblock Finding structure with randomness: Probabilistic algorithms for
  constructing approximate matrix decompositions.
\newblock {\em SIAM Review}, 53(2):217--288, 2011.

\bibitem{kahan1966numerical}
W.~Kahan.
\newblock Numerical linear algebra.
\newblock {\em Canadian Math. Bull}, 9(6):757--801, 1966.

\bibitem{martinsson2017randutv}
P.-G. Martinsson, G.~Quintana-Orti, and N.~Heavner.
\newblock rand{UTV}: A blocked randomized algorithm for computing a
  rank-revealing {UTV} factorization.
\newblock {\em arXiv preprint arXiv:1703.00998}, 2017.

\bibitem{stewart1992updating}
G.~W. Stewart.
\newblock An updating algorithm for subspace tracking.
\newblock {\em IEEE Transactions on Signal Processing}, 40(6):1535--1541, 1992.

\bibitem{stewart1993updating}
G.~W. Stewart.
\newblock Updating a rank-revealing {ULV} decomposition.
\newblock {\em SIAM Journal on Matrix Analysis and Applications},
  14(2):494--499, 1993.

\bibitem{stewart1999qlp}
G.~W. Stewart.
\newblock The {QLP} approximation to the singular value decomposition.
\newblock {\em SIAM Journal on Scientific Computing}, 20(4):1336--1348, 1999.

\end{thebibliography}

\end{document}